\documentclass[final]{amsart}
\usepackage{pdfsync}
\usepackage{amsmath,amssymb}
\usepackage{enumerate}
\usepackage{xspace}
\usepackage{boxedminipage}
\usepackage{algorithmicx}
\usepackage[ruled]{algorithm}
\usepackage{algpseudocode}
\usepackage{listings}
\usepackage{array}
\usepackage{tabularx}
\usepackage{subfig}
\usepackage{caption}
\usepackage{fullpage}
\usepackage{tristan}

\renewcommand{\bi}[2]{\ensuremath{\cA\qp{#1,#2}}}
\renewcommand{\bih}[2]{\ensuremath{\cA_h\qp{#1,#2}}}
\newcommand{\semih}[3]{\ensuremath{\cB\qp{#1;#2,#3}}}
\newcommand{\semi}[3]{\ensuremath{\cB\qp{#1;#2,#3}}}

\newcommand{\qnorm}[2]{\Norm{#1}_{\qp{#2}}}

\renewcommand{\vec}[1]{\geovec{#1}}
\renewcommand{\mat}[1]{\geomat{#1}}

\renewcommand{\enorm}[1]{\ensuremath{\norm{#1}}_{dG}}

\renewcommand{\sob}[2]{\ensuremath{\WW^{#1,#2}}}

\renewcommand{\leb}[1]{\ensuremath{\LL^{#1}}}
\newcommand{\dx}{\d \vec x}

\renewcommand{\fes}{V_h^k}

\numberwithin{equation}{section}

\setboolean{showtodo}{true}
\setboolean{shownotes}{true}
\setboolean{showchanges}{false}
\setboolean{usemathrsfs}{false}
\usepackage{tkz-graph}
\usetikzlibrary{arrows}

\author{
  James Jackaman
}
\address{
  James Jackaman
  \thanks{
    Department of Mathematics and Statistics, Whiteknights, University of Reading, Reading RG6 6AX, UK
    {\tt{James.Jackaman@pgr.reading.ac.uk}}.
}}

\author{
 Tristan Pryer
}
\address{
  Tristan Pryer
  \thanks{
    Department of Mathematics and Statistics, Whiteknights, University of Reading, Reading RG6 6AX, UK
    {\tt{T.Pryer@reading.ac.uk}}.
}}

\thanks{J.J. was supported through a PhD scholarship awarded by the
  ``EPSRC Centre for Doctoral Training in the Mathematics of Planet
  Earth at Imperial College London and the University of Reading''
  EP/L016613/1. T.P. was partially supported through the EPSRC grant
  EP/P000835/1 and the Newton Fund grant 261865400.}
\title[Quasinorms in semilinear elliptic problems]
      {Quasinorms in semilinear elliptic problems}
      \date{\today}

\begin{document}
\maketitle

\begin{abstract}
  In this note we examine the a priori and a posteriori analysis of
  discontinuous Galerkin finite element discretisations of semilinear elliptic PDEs with
  polynomial nonlinearity. We show that optimal a priori error bounds
  in the energy norm are only possible for low order elements using
  classical a priori error analysis techniques. We make use of
  appropriate quasinorms that results in optimal energy norm error
  control.

  We show that, contrary to the a priori case, a standard a posteriori
  analysis yields optimal upper bounds and does not require the
  introduction of quasinorms. We also summarise extensive numerical
  experiments verifying the analysis presented and examining the
  appearance of layers in the solution.
\end{abstract}

\section{Introduction}
\label{sec:introduction}

Let $\W \subset \reals^d$ with $d\geq 1$ be an open Lipschitz domain
and consider the problem: find $u\in \hoz$, such that
\begin{equation}
  \label{eq:pde}
  \begin{split}
    -\Delta u + \norm{u}^{p-2} u &= f \text{ in } \W
    \\
    u &= 0 \text{ on } \partial\W.
  \end{split}
\end{equation}
This class of equation is sometimes referred to as the
Lane-Emden-Fowler equation and are related to problems with critical
exponents \cite{clement1996quasilinear}. Furthermore, they arise in
the theory of boundary layers of viscous fluids
\cite{wong1975generalized}.

We are particularly interested in the class of PDE \eqref{eq:pde}
because of its application to the analysis of numerical schemes posed
for the KdV-like equation
\begin{equation}
  u_t
  -
  \qp{\norm{u}^{p-2} u}_x
  +
  u_{xxx} = 0.
  \label{eq:mkdv}
\end{equation}
Indeed, solutions of (\ref{eq:mkdv}) posed over a $1$-dimensional
domain satisfy
\begin{equation}
  \label{eq:energy}
    0
    =
    \ddt \cJ[u],
    \text{ with }      
  \cJ[u] = \int_{\W} \frac 1 2 \norm{u_x}^2 + \frac 1 p\norm{u}^p  \d x
\end{equation}
and energy minimisers of (\ref{eq:energy}) satisfy (\ref{eq:pde}) with
$f=0$ and appropriate boundary conditions. In
\cite{JackamanPapamikosPryer:2017,JackamanPryer:2018} a conservative
Galerkin scheme was proposed for (\ref{eq:mkdv}) and the a priori and
a posteriori analysis of this scheme requires quasi-optimal
approximation of the finite element solution of (\ref{eq:pde}) and
optimal a posteriori estimates. Hence our goal in this work is the
derivation and a priori and a posteriori bounds of Galerkin
discretisations of (\ref{eq:pde}).

We proceed as follows: In \S\ref{sec:setup} we introduce notation and
the model problem. We give some insight as to its properties that we
use in subsequent sections and propose a discontinuous Galerkin finite
element approximation. In \S\ref{sec:classical} we give a classical a
priori analysis based on arguments in \cite{Ciarlet:1978}. We show
that in the energy norm, the analysis is suboptimal for high
polynomial degrees and large values of $p$. In \S\ref{sec:quasinorms}
we modify the notion of a quasinorm from the works of
\cite{LiuBarrett:1996} to enable an optimal a priori error estimate to
be shown. In \S\ref{sec:apost} we derive an a posteriori estimate and
finally, in \S \ref{sec:numerics}, we showcase some numerical
experiments.

\section{Problem setup}
\label{sec:setup}

In this section we formulate the model problem, fix notation and give
some basic assumptions. Weakly, we may consider the PDE (\ref{eq:pde})
as: find $u\in \hoz$, such that
\begin{equation}
  \label{eq:weakform}
  \bi{u}{v} + \semi{u}{u}{v} = \ltwop{f}{v} \Foreach v\in\hoz,
\end{equation}
where $\ltwop{\cdot}{\cdot}$ denotes the $\leb{2}$ inner product and
the bilinear form $\cA : \hoz \times \hoz \to \reals$ is given by
\begin{equation}
  \label{eq:bilinear-form}
  \bi{u}{v} := \int_\W \nabla u\cdot\nabla v \dx.
\end{equation}
The semilinear form $\cB$ is given by
\begin{equation}
  \label{eq:semilinear-form}
  \semi{w}{u}{v} := \int_\W \norm{w}^{p-2} u v \dx.
\end{equation}
It is straightforward to verify this problem admits a unique solution.

\begin{Pro}[A priori bound 1]
  \label{eq:pro-apriorih1}
  Let $f\in\sobh{-1}(\W)$ and $u\in\hoz$ solve (\ref{eq:weakform}).
  Then we have
  \begin{equation}
    \frac 1 2 \Norm{\nabla u}_{\leb{2}(\W)}^2
    +
    \Norm{u}_{\leb{p}(\W)}^p
    \leq
    \frac 12 \Norm{f}_{\sobh{-1}(\W)}^2.    
  \end{equation}
\end{Pro}
\begin{Proof}[.]
  Using a standard energy argument, take $v=u$ in (\ref{eq:weakform}),
  then
  \begin{equation}
    \Norm{\nabla u}_{\leb{2}(\W)}^2
    +
    \Norm{u}_{\leb{p}(\W)}^p
    =
    \ltwop{f}{u}
    \leq
    \Norm{f}_{\sobh{-1}(\W)}\Norm{\nabla u}_{\leb{2}(\W)}
    \leq
    \frac 1 2 
    \qp{\Norm{f}_{\sobh{-1}(\W)}^2
      +
      \Norm{\nabla u}_{\leb{2}(\W)}^2},
  \end{equation}
  as required.
\end{Proof}

\begin{Pro}[A priori bound 2]
  \label{eq:pro-apriorilq}
  Let $f\in\leb{q}(\W)$, where $q = \frac{p}{p-1}$ and $u\in\hoz$
  solve (\ref{eq:weakform}) then we have
  \begin{equation}
    \Norm{\nabla u}_{\leb{2}(\W)}^2
    +
    \frac 1 q \Norm{u}_{\leb{p}(\W)}^p
    \leq
    \frac 1 q
    \Norm{f}_{\leb{q}(\W)}^q    
  \end{equation}
\end{Pro}
\begin{Proof}
  Again, take $v=u$ in (\ref{eq:weakform}), then
  \begin{equation}
    \Norm{\nabla u}_{\leb{2}(\W)}^2
    +
    \Norm{u}_{\leb{p}(\W)}^p
    =
    \ltwop{f}{u}
    \leq
    \Norm{f}_{\leb{q}(\W)}\Norm{u}_{\leb{p}(\W)}
    \leq
    \frac 1 p
    \Norm{u}_{\leb{p}(\W)}^p
    +
    \frac 1 q
    \Norm{f}_{\leb{q}(\W)}^q,
  \end{equation}
  as required.
\end{Proof}

\begin{Rem}[Behaviour of the a priori bounds in $p$]
  Note that the bounds given in Propositions \ref{eq:pro-apriorih1}
  and \ref{eq:pro-apriorilq} behave the same as $p$ increases, since
  \begin{equation}
    \lim_{p\to \infty}
    \frac 1 q
    =
    1,
  \end{equation}
  however the bound given in Proposition \ref{eq:pro-apriorilq} blows
  up as $p$ decreases, indeed 
  \begin{equation}
    \lim_{p\to 1}
    \frac 1 q
    =
    \infty.
  \end{equation}
  We will only consider the case $p \geq 2$ in this work.
\end{Rem}

\subsection{Discretisation}

Let $\T{}$ be a regular subdivision of $\Omega$ into disjoint
simplicial elements. We assume that the subdivision $\T{}$ is
shape-regular \cite[p.124]{Ciarlet:1978}, is
$\closure{\Omega}=\cup_{K}\closure{K}$ and that the elemental faces
are points (for $d=1$), straight lines (for $d=2$) or planar (for
$d=3$) segments; these will be, henceforth, referred to as
\emph{facets}.  By $\Gamma$ we shall denote the union of all
($d-1$)-dimensional facets associated with the subdivision $\T{}$
including the boundary. Further, we set
$\Gamma_{i}:=\Gamma\backslash\partial\Omega$.

For a nonnegative integer $k$, we denote the set of all polynomials of
total degree at most $k$ by $\poly{k}(K)$. For $k \geq 1$, we consider the
finite element space
\begin{equation}
 \fes :=\{v\in L^2(\Omega):v|_{K} \in\poly{k}(K), \forall K\}.
\end{equation}

Further, let $K^+$, $K^-$ be two (generic) elements sharing a facet
$e:=\partial K^+\cap\partial K^-\subset\Gamma_{i}$ with respective
outward normal unit vectors $\vec{n}^+$ and $\vec{n}^-$ on $e$. For a
function $v:\Omega\to\reals$ that may be discontinuous across
$\Gamma_{i}$, we set $v^+:=v|_{e\subset\partial K^+}$,
$v^-:=v|_{e\subset\partial K^-}$, and we define the jump by
\[ {\jump{v}:=v^+\vec{n}^++v^-\vec{n}^-};\]
if $e\in \partial K\cap\partial\Omega$, we set
$\jump{v}:=v^+\vec{n}$. Also, we define $h_{K}:=\diam(K)$ and we
collect them into the element-wise constant function $
h:\Omega\to\mathbb{R}$, with $h|_{K}=h_{K}$, $K$, $
h|_e=({h_{K^+}+h_{K^-}})/2$ for $e\subset\Gamma_{i}$ and $h |_e=h_K$
for $e\subset\partial K\cap\partial\Omega$. We assume that the
families of meshes considered in this work are locally
quasi-uniform. Note that this restriction can be relaxed by following
arguments as in \cite{GeorgoulisMakridakisPryer:2017}.

For $s>0$, we define the \emph{broken} Sobolev space $\sobh{s}(\T{})$, by
\[
\sobh{s}(\T{}):= \{w\in \leb{2}(\W): w|_K\in \sobh{s}(K), K\in\T{}\},
\]
along with the broken (element-wise) gradient and Laplacian $\nabla_h\equiv
\nabla_h(\T{})$ and $\Delta_h\equiv \Delta_h(\T{})$.

We consider the interior penalty (IP) discontinuous Galerkin
discretisation of (\ref{eq:bilinear-form}), reading: find $u_h \in
\fes$ such that
\begin{equation}\label{eq:dg}
  \bih{u_h}{v_h}
  +
  \semih{u_h}{u_h}{v_h} = \ltwop{f}{v_h} \Foreach v_h \in \fes,
\end{equation}
where  
\begin{equation}
  \label{eq:IP}
  \begin{split}
    \bih{u_h}{v_h} 
    &=
    \int_\W \nabla_h u_h \cdot \nabla_h v_h \dx
    -
    \int_{\Gamma} \big(\jump{v_h} \cdot \avg{P_{k-1}\qp{\nabla u_h}}
    +
    \jump{u_h} \cdot \avg{P_{k-1}\qp{\nabla v_h}}
    -
    \sigma \jump{u_h}\cdot \jump{v_h}\big)\d s
    ,
  \end{split}
\end{equation}
where $\sigma > 0$ is the, so-called, \emph{discontinuity penalisation
  parameter} given by
\begin{equation}
  \label{eq:sigma}
  \sigma:= C_{\sigma}\frac{k^2}{h}
\end{equation}
and $P_{k-1} : \leb{2}(\W) \to V_h^{k-1}$ denotes the $\leb{2}$
orthogonal projection operator. This is included in the bilinear form
to ensure that $\cA_h$ is well defined over $\sobh{1}(\W)\times
\sobh{1}(\W)$.

\begin{Defn}[Mesh dependent norms]
  \label{def:mesh-dep-norms}
  We introduce the {mesh dependent} $\sobh1$ norm to be
  \begin{gather}
    \enorm{w}^2 := \Norm{\nabla_h w}_{\leb{2}(\W)}^2 + \Norm{\sqrt{\sigma}\jump{w}}_{\leb{2}(\Gamma)}^2.
  \end{gather} 
\end{Defn}

Note that the the bilinear form (\ref{eq:bilinear-form}) satisfies
boundedness and coercivity properties for $C_\sigma$ chosen large
enough \cite[c.f.]{ErnGuermond:2004}, that is
\begin{equation}
  \label{eq:dgbdd}
  \begin{split}
    \bih{u}{v} &\leq \widetilde C_B \enorm{u}\enorm{v}
    \\
    \bih{u}{u} &\geq \widetilde C_C \enorm{u}^2 \Foreach u, v \in\fes + \sobh1_0(\W).
  \end{split}
\end{equation}

\section{Classical a priori analysis}
\label{sec:classical}

In this section we examine analysis based on classical arguments such
as those used in \cite{Ciarlet:1978} for the $p$-Laplacian.

\begin{Lem}[{Properties of $\cB(\cdot; \cdot, \cdot)$, cf.\ \cite[\S5.3]{Ciarlet:1978}}]
  \label{lem:sandri}
  There exist constants
  \begin{enumerate}
  \item  $C_L>0$ such that
    \begin{equation}
      \semi{u-u_h}{u-u_h}{u-u_h}
      \leq
      C_L    
      \qp{\semi{u}{u}{u-u_h}
      -
      \semi{u_h}{u_h}{u-u_h}}
    \end{equation}
  \item $C_U>0$ such that
    \begin{equation}
      \semi{u}{u}{u-v_h}
      -
      \semi{u_h}{u_h}{u-v_h}
      \leq
      C_U
      \Norm{u-u_h}_{\leb{p}(\W)} 
      \Norm{u-v_h}_{\leb{q}(\W)} 
      .
    \end{equation}
  \end{enumerate}
\end{Lem}

\begin{The}
  \label{the:best-approx}
  Let $u\in\sobh{2}(\W)\cap\hoz$ solve (\ref{eq:pde}) and $u_h
  \in\fes$ be the finite element approximation of (\ref{eq:dg}) then
  for $k\geq 1$ we have
  \begin{equation}
    \enorm{u - u_h}^2 + \Norm{u - u_h}_{\leb{p}(\W)}^p
    \leq
    C
    \inf_{v_h\in\fes\cap \cont{0}(\W)}
    \qp{
      \enorm{u - v_h}^2 + \Norm{u - v_h}_{\leb{q}(\W)}^q
    },
  \end{equation}
  where $q = \frac{p}{p-1}$ is the Sobolev conjugate of $p$.
\end{The}
\begin{Proof}[.]
  Since $u_h$ solves (\ref{eq:dg}), we have, though Lemma
  \ref{lem:sandri} and Galerkin orthogonality over $\fes\cap \cont{0}(\W)$
  \begin{equation}
    \label{eq:appf1}
    \begin{split}
      \widetilde C_C
      \enorm{u - u_h}^2
      +
      \frac1{C_L}
      \Norm{u - u_h}_{\leb{p}(\W)}^p
      &\leq
      \bih{u - u_h}{u - u_h}
      +
      \frac1{C_L}
      \semih{u - u_h}{u - u_h}{u - u_h}
      \\
      &\leq
      \bih{u - u_h}{u - u_h}
      +
      \semih{u}{u}{u-u_h}
      - 
      \semih{u_h}{u_h}{u-u_h}
      \\
      &\leq
      \bih{u - u_h}{u - v_h}
      +
      \semih{u}{u}{u-v_h}
      - 
      \semih{u_h}{u_h}{u-v_h},
    \end{split}
  \end{equation}
  for any $v_h \in \fes\cap \cont{0}(\W)$.
  Note that
  \begin{equation}
    \label{eq:appf2}
    \begin{split}
      \bih{u - u_h}{u - v_h}
      &\leq
      \frac {\widetilde C_B^2}{2 \widetilde C_C }
      \enorm{u - v_h}^2
      +
      \frac {\widetilde C_C}  2
      \enorm{u - u_h}^2.
    \end{split}
  \end{equation}
  Further,
  \begin{equation}
    \label{eq:appf3}
    \begin{split}
      \semih{u}{u}{u-v_h}
      - 
      \semih{u_h}{u_h}{u-v_h}
      &\leq
      C_U
      \Norm{u - u_h}_{\leb{p}(\W)}
      \Norm{u - v_h}_{\leb{q}(\W)}
      \\
      &\leq
      \frac 1{2C_L}
      \Norm{u - u_h}_{\leb{p}(\W)}^p
      +
      \qp{\frac p {2C_L}}^{-q/p}\frac{C_U}{q}
      \Norm{u - v_h}_{\leb{q}(\W)}^q.
    \end{split}
  \end{equation}
  Substituting (\ref{eq:appf2}), (\ref{eq:appf3}) into (\ref{eq:appf1})
  yields the desired result.
\end{Proof}

\begin{Cor}
  \label{eq:cor-energy}
  Choosing $v_h = I_k u$, the Cl\'ement interpolant of $u$, in Theorem
  \ref{the:best-approx} and under further smoothness requirements,
  that $u\in\sob{k+1}{p}(\W)$, we see that
  \begin{equation}
    \enorm{u - u_h}^2 + \Norm{u - u_h}_{\leb{p}(\W)}^p
    \leq
    C
    \qp{h^{2k}\norm{u}_{\sobh{k+1}(\W)} + h^{\qp{k+1}q}\norm{u}_{\sob{k+1}{p}(\W)}}.
  \end{equation}
\end{Cor}

\begin{Rem}[Optimality of Corollary \ref{eq:cor-energy}]
  Notice that the bound given in Corollary \ref{eq:cor-energy} depends
  upon $q = \frac{p}{p-1}$. Notice, as shown in Table
  \ref{tab:table-opt-energy}, the energy error bounds are optimal only
  if $p=2$ for all $k$ or $k=1$ or all $p$.
  \begin{table}[h!]
  \begin{center}
    \caption{ \label{tab:table-opt-energy} In the following table we
      examine the optimality of the finite element approximation in
      the energy norm. The numerical values in the table correspond to
      $\min\qp{k,\frac{\qp{k+1}q}2}$ and are coloured green or red depending upon
      whether the bound is optimal, in the function approximation
      sense, or suboptimal respectively.}
    \begin{tabular}{|c|c|c|c|c|c|}
      \hline
      $k$&  $p=2$ &  $p=3$ &  $p=4$ &  $p=5$ &  $p\to\infty$
      \\
      \hline
      $1$ & \textcolor{green}{$2$} & \textcolor{green}{$3/2$} & \textcolor{green}{$4/3$} & \textcolor{green}{$5/4$} & \textcolor{green}{$1$}
      \\
      $2$ & \textcolor{green}{$3$} & \textcolor{green}{$9/4$} & \textcolor{green}{$2$} & \textcolor{red}{${15}/8$} & \textcolor{red}{$3/2$}
      \\
      $3$ & \textcolor{green}{$4$} & \textcolor{green}{$3$} & \textcolor{red}{$8/3$} & \textcolor{red}{$5/2$} & \textcolor{red}{$2$}
      \\
      $4$ & \textcolor{green}{$5$} & \textcolor{red}{${15}/4$} & \textcolor{red}{${10}/3$} & \textcolor{red}{${25}/8$} & \textcolor{red}{$5/2$}
      \\
      \hline
    \end{tabular}
  \end{center}
\end{table}
\end{Rem}

\begin{Rem}[Dual bounds]
  This lack of optimality propagates further when consider bounds
  based on duality approaches. Indeed, using the dual problem
  \begin{equation}
    \begin{split}
      -\Delta z + \qp{p-1} u^{p-2} z &= u - u_h \text{ in } \W
      \\
      z &= 0 \text{ on }\partial \W,
    \end{split}
  \end{equation}
  one can show that 
  \begin{equation}
    \Norm{u - u_h}_{\leb{2}(\W)}
    \leq
    C
    \qp{
      h\enorm{u-u_h}
      +
      \Norm{u - u_h}_{\leb{p}(\W)}^2
    }.
  \end{equation}
  We will not prove this here for brevity but, as illustrated in Table
  \ref{tab:table-opt-dual}, the bound is optimal only when $p=2$.
\end{Rem}



  \begin{table}[h!]
  \begin{center}
    \caption{ \label{tab:table-opt-dual} In the following table we
      examine the optimality of the finite element approximation in
      the $\leb{2}$ norm. The numerical values in the table correspond
      to $\min\qp{{k+1}, {\qp{k+1}q/2 + 1}, {4k/p}, {(2k+2)/(p-1)}}$
      thus represent the convergence rate in $\leb{2}$. They are
      coloured green or red depending upon whether the bound is
      optimal, in the function approximation sense, or suboptimal
      respectively. Notice that they $\leb{2}$ norm estimate is only
      optimal for $p=2$ for all $k$.}
    \begin{tabular}{|c|c|c|c|c|c|}
      \hline
      $k$&  $p=2$ &  $p=3$ &  $p=4$ &  $p=5$ &  $p\to\infty$
      \\
      \hline
      $1$ & \textcolor{green}{$2$} & \textcolor{red}{$4/3$} & \textcolor{red}{$1$} & \textcolor{red}{$4/5$} & \textcolor{red}{$0$}
      \\
      $2$ & \textcolor{green}{$3$} & \textcolor{red}{$8/3$} & \textcolor{red}{$2$} & \textcolor{red}{$3/2$} & \textcolor{red}{$0$}
      \\
      $3$ & \textcolor{green}{$4$} & \textcolor{green}{$4$} & \textcolor{red}{$8/3$} & \textcolor{red}{$2$} & \textcolor{red}{$0$}
      \\
      $4$ & \textcolor{green}{$5$} & \textcolor{red}{${19}/4$} & \textcolor{red}{${10}/3$} & \textcolor{red}{${5}/2$} & \textcolor{red}{$0$}
      \\
      \hline
    \end{tabular}
  \end{center}
\end{table}

\section{A priori analysis based on quasi norms}
\label{sec:quasinorms}

In this section we will examine the use of quasinorms to rectify the
gap in the a priori analysis.

\begin{Defn}[Quasinorm]
  Let $v\in\leb{p}(\W)$, $p\geq 2$, then for any $w\in\leb{p}(\W)$ we
  define the quasinorm
  \begin{equation}
    \qnorm{v}{w,p}^2 := \int_\W \norm{v}^2 \qp{\norm{w} + \norm{v}}^{p-2} \d x.
  \end{equation}
  This satisfies the usual properties of a norm, in that
  \begin{equation}
      \qnorm{v}{w,p} \geq 0 \text{ and } \qnorm{v}{w,p} = 0 \iff v = 0.
  \end{equation}
  However, the usual triangle inequality is replaced by
  \begin{equation}
    \qnorm{v_1 + v_2}{w,p} \leq C\qp{\qnorm{v_1}{w,p}
      +
      \qnorm{v_2}{w,p}
    },
  \end{equation}
  where $C = C(v_1, v_2, w, p)$
\end{Defn}

\begin{Rem}[Properties of the quasinorm]
  \label{rem:properties-quasinorm}
  As can be seen from the definition, the quasinorm is related to the
  $\leb{p}$ norm through the relationship
  \begin{equation}
    \Norm{v}_{\leb{p}(\W)}^p
    \leq
    \qnorm{v}{w,p}^2
    \leq
    C
    \Norm{v}_{\leb{p}(\W)}^2,
  \end{equation}
  for $v\in\leb{p}(\W)$, $p\geq 2$ and any $w\in\leb{p}(\W)$. The key
  property that the quasinorm satisfies that allows for optimal a
  priori treatment is that the semilinear form is coercive with
  respect to it, that is
  \begin{equation}
    \label{eq:coer}
    \semi{u}{u}{u-v}
    -
    \semi{v}{v}{u-v}
    \geq
    \overline C_C
    \qnorm{u-v}{u,p}^2.
  \end{equation}
  In addition, it is bounded \cite{ebmeyer2005quasi} in that for any
  $\theta>0$ there exists a $\gamma>0$ such that
  \begin{equation}
    \label{eq:bdd}
    \norm{\semi{u}{u}{w}
      -
      \semi{v}{v}{w}}
    \leq
    \overline C_B
    \qp{\theta^\gamma \qnorm{u-v}{u,p}^2
      +
      \theta \qnorm{w}{u,p}^2
    },
  \end{equation}
  where
  \begin{equation}
    \gamma =
    \begin{cases}
      1 \text{ if } \theta < 1
      \\
      \frac 1 {p-1}
      \text{ if } \theta \geq 1.
    \end{cases}
  \end{equation}
  It was the lack of a sufficiently sharp boundedness property that
  led to suboptimality in the analysis presented in Section
  \ref{sec:classical}. The key observation to rectify the
  suboptimality is to measure the error in the norm
  \begin{equation}
    \enorm{u - u_h}^2 + \qnorm{u-u_h}{u,p}^2,
  \end{equation}
  rather than the energy norm.

  Henceforth, we will use the notation
  \begin{equation}
    C_C := \min\qp{\overline C_C, \widetilde C_C}
    \qquad
    C_B := \max\qp{\overline C_B , \widetilde C_B}.
  \end{equation}
\end{Rem}

\begin{Pro}[A priori bound 3]
  Let $f\in\sobh{-1}(\W)$ and $u\in\hoz$ solve (\ref{eq:weakform})
  then
  \begin{equation}
    \Norm{\nabla u}^2_{\leb{2}(\W)}
    +
    2^{3-p}\qnorm{u}{u,p}^2
    \leq
    \Norm{f}_{\sobh{-1}(\W)}^2.
  \end{equation}
\end{Pro}
\begin{Proof}
  Notice that
  \begin{equation}
    \Norm{\nabla u}^2_{\leb{2}(\W)}
    +
    2^{2-p}\qnorm{u}{u,p}^2
    =
    \bi{u}{u} + \semi{u}{u}{u}
    =
    \ltwop{f}{u}
    \leq
    \frac 1 2 \qp{
      \Norm{f}_{\sobh{-1}(\W)}^2 +
    \Norm{\nabla u}_{\leb{2}(\W)}},
  \end{equation}
  as required.
\end{Proof}

\begin{The}
  \label{the:best-approx-quasi}
  Let $u\in\sobh{2}(\W)\cap\hoz$ solve (\ref{eq:pde}) and $u_h$ be the finite element
  approximation of (\ref{eq:dg}) then for $k\geq 1$ we have
  \begin{equation}
    \enorm{u - u_h}^2 + \qnorm{u - u_h}{u,p}^2
    \leq
    C
    \inf_{v_h\in\fes}
    \qp{
      \enorm{u - v_h}^2 + \qnorm{u - v_h}{u,p}^2
    }.
  \end{equation}
\end{The}

\begin{Proof}
  Making use of the coercivity of $\cB$ we have
  \begin{equation}
    \begin{split}
      C_C\qp{
        \enorm{u - u_h}^2
        +
        \qnorm{u - u_h}{u,p}^2
      }
      &\leq
      \bi{u - u_h}{u - u_h}
      +
      \semi{u}{u}{u-u_h}
      -
      \semi{u_h}{u_h}{u-u_h}
      \\
      &=
      \bi{u - u_h}{u - v_h}
      +
      \semi{u}{u}{u-v_h}
      -
      \semi{u_h}{u_h}{u-v_h},
    \end{split}
  \end{equation}
  for any $v_h\in\fes\cap\cont{0}(\W)$, using Galerkin orthogonality. Now, through
  (\ref{eq:bdd}) and (\ref{eq:dgbdd}) we have
  \begin{equation}
    \begin{split}
      C_C
      \qp{
        \enorm{u - u_h}^2
        +
        \qnorm{u - u_h}{u,p}^2
      }
      &\leq
      \frac{C_C}2 \enorm{u - u_h}^2
      +
      \frac{C_B^2}{2C_C} \enorm{u - v_h}^2
      \\
      &\qquad +
      C_B
      \qp{\theta^\gamma \qnorm{u-u_h}{u,p}^2
        +
        \theta \qnorm{u-v_h}{u,p}^2
      }.
    \end{split}
  \end{equation}
  Choosing $\theta = \min\qp{\frac{C_C}{2C_B}, \frac 1 2}$ then $\gamma =
    1$. Rearranging the inequality yields the desired result.
\end{Proof}

\begin{Lem}
  \label{lem:approximation-quasi}
  Let $p \geq 2$ and $v \in \sob{k+1}{p}(\W)$ then
  \begin{equation}
    \inf_{v_h\in\fes} \qnorm{v-v_h}{v,p}
    \leq
    C h^{k+1} \norm{v}_{\sob{k+1}{p}(\W)}.
  \end{equation}
\end{Lem}
\begin{Proof}
  Using the property of the quasinorm given in Remark
  \ref{rem:properties-quasinorm} we have
  \begin{equation}
    \qnorm{v-v_h}{v,p}^2
    \leq
    C
    \Norm{v-v_h}_{\leb{p}(\W)}^2,
  \end{equation}
  and the result follows from best approximation in $\leb{p}(\W)$.
\end{Proof}

\begin{Cor}
  \label{cor:quasi}
  Under the conditions of Theorem \ref{the:best-approx-quasi} suppose
  that $u\in \sobh{k+1}(\W)\cap \hoz \cap \sob{k}{p}(\W)$, then
  \begin{equation}
    \qp{\enorm{u - u_h}^2 + \qnorm{u - u_h}{u,p}^2}^{1/2}
    \leq
    C
    h^k
    \qp{\norm{u}_{\sobh{k+1}(\W)}
      + \norm{u}_{\sob{k}{p}(\W)}
    }.
  \end{equation}
\end{Cor}

\begin{Rem}[Optimality of Corollary \ref{cor:quasi}]
  Notice that the bound given in Corollary \ref{cor:quasi} is optimal
  regardless of the choice of $p$ for smooth enough $u$.
\end{Rem}

\begin{Rem}[Dual bounds]
  By modifying the dual problem to
  \begin{equation}
    \begin{split}
      -\Delta z + \qp{p-1} u^{p-2} z &= \qp{u - u_h}\qp{\norm{u} + \norm{u-u_h}}^{p-2} \text{ in } \W
      \\
      z &= 0 \text{ on }\partial \W,
    \end{split}
  \end{equation}
  one can also show optimal a priori bounds for the quasinorm error.
\end{Rem}

\section{A posteriori error analysis}
\label{sec:apost}

In this section we derive a reliable a posteriori estimator.

\begin{Pro}[A priori bound 4]
  \label{pro:anotherbound}
  Let $f\in\sobh{-1}(\W)$ and $u\in\hoz$ solve \eqref{eq:weakform} and $w\in\hoz$
  solve
  \begin{equation}
    \label{eq:pert}
    \bi{w}{v} + \semi{w}{w}{v} = \ltwop{f - \mathfrak R}{v} \Foreach v\in\hoz,
  \end{equation}
  for some $\mathfrak R \in \sobh{-1}(\W)$ then
  \begin{equation}
    \Norm{\nabla u - \nabla w}_{\leb{2}(\W)}^2
    +
    2C_L \Norm{u-w}_{\leb{p}(\W)}^p
    \leq
    \Norm{\mathfrak R}_{\sobh{-1}(\W)}^2
  \end{equation}
\end{Pro}
\begin{Proof}
  Through the definitions of $u$ and $w$, we have the relation that
  \begin{equation}
    \bi{u - w}{v} + \semi{u}{u}{v} - \semi{w}{w}{v} = \ltwop{\mathfrak R}{v} \Foreach v\in\hoz.
  \end{equation}
  Hence, choosing $v=u-w$
  \begin{equation}
    \Norm{\nabla u - \nabla w}^2_{\leb{2}(\W)} + C_L\Norm{u - w}_{\leb{p}(\W)}^p
    \leq
    \ltwop{\mathfrak R}{u-w}
    \leq
    \frac1 2 \qp{    \Norm{\mathfrak R}_{\sobh{-1}(\W)}^2
      +
      \Norm{\nabla u - \nabla w}^2_{\leb{2}(\W)}},
  \end{equation}
  as required.
\end{Proof}

To invoke the results of Proposition \ref{pro:anotherbound} we require
an object $w\in\sobh{1}(\W)$. The dG solution $u_h\notin \sobh{1}(\W)$, so we make
use of an appropriate postprocessor as an intermediate quantity.

\begin{Lem}[\cite{KarakashianPascal:2003}]
  \label{lem:KP}
  Let $\mathcal{N}$ denote the set of all Lagrange nodes of $\fes$,
  and $\E :\fes\to \fes\cap H^1_0(\Omega)$ be defined on the
  conforming Lagrange nodes $\nu\in\mathcal{N}$ by
  \[
  \E(v)(\nu):= \bigg \{ \begin{array}{cc}
    \displaystyle|\omega_{\nu}|^{-1}\sum_{K\in\omega_{\nu}} v|_{K}({\nu}),& \nu\in\Omega; \\ 
    0, & \nu\in \partial\Omega,
  \end{array} 
  \]
  with $ \omega_{\nu}:=\bigcup_{K\in\T{}: \nu\in\overline{K}}K, $ the
  set of elements sharing the node $\nu\in\mathcal{N}$ and
  $|\omega_{\nu}|$ their cardinality. Then, the following bound holds
  \begin{equation}\label{KP_stab}
    \sum_{K\in\T{}}
    \norm{v - \E(v)}_{\sobh{\alpha}(K)}^2
    \leq
    C_{\alpha}
    \sum_{e\in\Gamma} \Norm{{h}^{1/2-\alpha}\jump{v}}_{\leb{2}(e)}^2,
  \end{equation}
  with $\alpha=0,1$, $C_{\alpha}\equiv C_{\alpha}(k)>0$ a
  constant independent of $h$, $v$ and $\T{}$, but depending on the
  shape-regularity of $\T{}$ and on the polynomial degree
  $k$. 
\end{Lem}

\begin{Pro}
  \label{pro:recon-pde}
  The reconstruction $\E(u_h)$ satisfies the perturbed PDE
  \begin{equation}
    \bi{\E(u_h)}{v} + \semi{\E(u_h)}{\E(u_h)}{v} = \ltwop{f - \mathfrak R}{v} \Foreach v\in\hoz,
  \end{equation}
  with
  \begin{equation}
    \begin{split}
      \ltwop{\mathfrak R}{v}
      &=
      \bih{u_h - \E(u_h)}{v} + \semi{u_h}{u_h}{v} - \semi{\E(u_h)}{\E(u_h)}{v}
      \\
      &\qquad +
      \ltwop{f}{v-v_h} - \bih{u_h}{v-v_h} - \semi{u_h}{u_h}{v-v_h}
      \Foreach v_h\in \fes.
    \end{split}
  \end{equation}
\end{Pro}

\begin{The}
  \label{the:apost-E}
  Let $f\in\leb{2}(\W)$ and $u\in\hoz$ solve \eqref{eq:weakform} and
  $u_h\in\fes$ solve (\ref{eq:dg}). Further let $\E(u_h)\in\hoz$
  denote the reconstruction operator given in Lemma
  \ref{lem:KP}. Then,
  \begin{equation}
    \Norm{\nabla u - \nabla \E(u_h)}_{\leb{2}(\W)}^2
    +
    2C_L \Norm{u-\E(u_h)}_{\leb{p}(\W)}^p
    \leq
    C\sum_{K\in\T{}} \qb{ \eta_R^2
    +
    \sum_{e\in \partial K} \eta_J^2},
  \end{equation}
  where
  \begin{equation}
    \begin{split}
      \eta_R^2
      &:=
      \Norm{h\qp{f+\Delta u_h - \norm{u_h}^{p-2} u_h}}_{\leb{2}(K)}^2
      \\
      \eta_J^2
      &:=
      \Norm{h^{1/2}\jump{\nabla u_h}}_{\leb{2}(e)}^2
      +
      \Norm{h^{-1/2}\jump{u_h}}_{\leb{2}(e)}^2
      +
      \Norm{h^{1/2}\jump{u_h}}_{\leb{2}(e)}^2
    \end{split}
  \end{equation}
\end{The}

\begin{Proof}[.]
  It suffices to determine an upper bound for $\Norm{\mathfrak
    R}_{\sobh{-1}(\W)}$. To that end, by Proposition
  \ref{pro:recon-pde}
  \begin{equation}
    \begin{split}
      \ltwop{\mathfrak R}{v}
      &=
      \underbrace{\bih{u_h - \E(u_h)}{v}}_{\cI_1}
      +
      \underbrace{\semi{u_h}{u_h}{v} - \semi{\E(u_h)}{\E(u_h)}{v}}_{\cI_2}
      \\
      &\qquad +
      \underbrace{\ltwop{f}{v-v_h} - \bih{u_h}{v-v_h} - \semi{u_h}{u_h}{v-v_h}}_{\cI_3},
    \end{split}
  \end{equation}
  and we proceed to bound the terms individually. Firstly,
  \begin{equation}
    \label{eq:first}
    \begin{split}
      \cI_1
      &=
      \sum_{K\in\T{}}
      \int_K \qp{\nabla u_h - \nabla \E(u_h)} \cdot \nabla v \dx
      -
      \sum_{e\in\Gamma}
      \jump{u_h} \avg{P_{k-1} \qp{\nabla v}} \d s
      \\
      &\leq
      \sum_{K\in\T{}}
      \Norm{\nabla u_h - \nabla\E(u_h)}_{\leb{2}(K)}
      \Norm{\nabla v}_{\leb{2}(K)}
      +
      \sum_{e\in\Gamma} \Norm{h^{-1/2} \jump{u_h}}_{\leb{2}(e)} \Norm{h^{1/2} \avg{P_{k-1} \qp{\nabla v}}}_{\leb{2}(e)}
      \\
      &\leq
      \qp{C_1^{1/2}+C_{dim,\T{}}^{1/2}C_{trace}^{1/2}}\qp{
        \sum_{e \in \Gamma}
        \Norm{h^{-1/2} \jump{u_h}}_{\leb{2}(e)}^2
      }^{1/2}
      \qp{
        \sum_{K\in\T{}}
        \Norm{\nabla v}_{\leb{2}(K)}^2
      }^{1/2}
      \\
      &\leq
      C\qp{
        \sum_{K\in\T{}}
        \sum_{e \in \partial K}
        \eta_J^2
      }^{1/2}
      \Norm{\nabla v}_{\leb{2}(\W)},
    \end{split}
  \end{equation}
  where $C_{dim,\T{}}^{1/2}$ is a constant depending on the dimension
  and the triangulation and $C_{trace}$ is the constant from a trace
  estimate. The second term can be controlled by
  \begin{equation}
    \label{eq:second}
    \begin{split}
      \cI_2
      &=
      \int_\W \qp{\norm{u_h}^{p-2} u_h - \norm{\E(u_h)}^{p-2} \E(u_h)} v \dx
      \\
      &\leq
      C(u_h, \E(u_h), p)
      \sum_{K\in\T{}} \Norm{u_h - \E(u_h)}_{\leb{2}(K)} \Norm{v}_{\leb{2}(K)}.
      \\
      &\leq
      C(u_h, \E(u_h), p)
      \qp{\sum_{K\in\T{}} \Norm{u_h - \E(u_h)}^2_{\leb{2}(K)}}^{1/2} \Norm{v}_{\leb{2}(\W)}
      \\
      &\leq
      C_P C(u_h, \E(u_h), p) C_0^{1/2}
      \qp{\sum_{e\in\Gamma} \Norm{h^{1/2} \jump{u_h}}^2_{\leb{2}(e)}}^{1/2} \Norm{\nabla v}_{\leb{2}(\W)},
    \end{split}
  \end{equation}
  where $C_P$ is the Poincar\'e constant. To finish, $\cI_3$ is
  controlled by a standard a posteriori argument.
  \begin{equation}
    \label{eq:3.24}
    \begin{split}
      \cI_3
      &=
      \int_\W f\qp{v - v_h} - \nabla_h u_h \cdot \qp{\nabla v - \nabla_h v_h} - \norm{u_h}^{p-2} u_h \qp{v-v_h}\d \vec x
      +
      \int_\Gamma \jump{u_h} \cdot \avg{P_{k-1}\qp{\nabla v - \nabla_h v_h}} \d s
      \\&\qquad 
      + 
      \int_\Gamma
      \jump{v -v_h} \avg{\nabla_h u_h}
      -
      \sigma \jump{u_h} \cdot \jump{{v - v_h}} \d s
      \\
      &=
      \int_\W \qp{f + \Delta_h u_h - \norm{u_h}^{p-2} u_h} \qp{v - v_h} \d \vec x
      -
      \int_\Gamma \jump{\nabla u_h} \avg{v - v_h} \d s
      \\
      & \qquad + 
      \int_\Gamma \jump{u_h} \cdot \avg{P_{k-1}\qp{\nabla v - \nabla_h v_h}}
      -
      \sigma \jump{u_h} \cdot \jump{{v - v_h}} \d s.
    \end{split}
  \end{equation}
  Splitting the integrals elementwise and making use of the
  Cauchy-Schwarz inequality we see
  \begin{equation}
    \label{eq:pf2}
    \begin{split}
      \int_\W \qp{f + \Delta u_h - \norm{u_h}^{p-2} u_h} \qp{v - v_h} \d \vec x
      &\leq
      \sum_{K\in\T{}}
      \Norm{h\qp{f + \Delta u_h - \norm{u_h}^{p-2} u_h}}_{\leb{2}(K)} 
      \Norm{h^{-1}\qp{v - v_h}}_{\leb{2}(K)} 
    \end{split}
  \end{equation}
  Similarly, for the second,
  \begin{equation}
    \label{eq:pf3}
    \begin{split}
      -\sum_{e\in\Gamma} \int_e \jump{\nabla u_h} \avg{v - v_h} \d s
      &\leq
      \sum_{K\in\T{}} \qb{
        \sum_{e\in \partial K}
        \Norm{h^{1/2}\jump{\nabla u_h}}_{\leb{2}(e)} \Norm{h^{-1/2} \avg{v - v_h}}_{\leb{2}(e)}
      },
    \end{split}
  \end{equation}
  and third term
  \begin{equation}
    \label{eq:pf4}
    \begin{split}
      \sum_{e\in\Gamma}
      \int_e \jump{u_h} \cdot \avg{{P_{k-1}\qp{\nabla v - \nabla v_h}}} \d s
      &\leq
      \sum_{K\in\T{}}
      \qb{
        \sum_{e\in \partial K}
        \Norm{h^{-1/2} \jump{u_h}}_{\leb{2}(e)}
        \Norm{h^{1/2} \avg{P_{k-1}\qp{\nabla v - \nabla v_h}}}_{\leb{2}(e)}
      }.
    \end{split}
  \end{equation}
  For the final term
  \begin{equation}
    \label{eq:pf5}
    \begin{split}
      \sum_{e\in\Gamma}
      \int_e \sigma \jump{u_h} \cdot \jump{{v - v_h}} \d s
      &\leq 
      C_\sigma
      \sum_{K\in\T{}}\qb{
        \sum_{e\in \partial K}
        \Norm{h^{-1/2}\jump{u_h}}_{\leb{2}(e)}
        \Norm{h^{-1/2}\jump{{v - v_h}}}_{\leb{2}(e)}
      }
    \end{split}
  \end{equation}
  Collecting (\ref{eq:3.24})--(\ref{eq:pf5}) we have
  \begin{equation}
    \begin{split}
      \cI_3
      &\leq
      C \qp{\sum_{K\in\T{}} \qb{\eta_{R} + \sum_{e\in\partial K}\eta_{J}}  \Phi(v-v_h)},
    \end{split}
  \end{equation}
  where
  \begin{equation}
    \begin{split}
      \Phi(w)
      &=
      \max\bigg(
      \Norm{h^{-1} {w}}_{\leb{2}(K)}
      ,
      \max_{e\in\Gamma} \Norm{h^{-1/2} {w}}_{\leb{2}(e)}
      ,
      \max_{e\in\Gamma} \Norm{h^{1/2} \avg{P_{k-1}\qp{\nabla w}}}_{\leb{2}(e)}
      ,
      \max_{e\in\Gamma} \Norm{h^{-1/2}\jump{w}}_{\leb{2}(e)}
      \bigg).
    \end{split}
  \end{equation}
  Choosing $v_h = P_0 v$ and in view of the approximation properties
  and stability of the $\leb{2}$ projector we have that
  \begin{equation}
    \Phi(v-v_h) \leq C \Norm{\nabla v}_{\leb{2}(\patch K)},
  \end{equation}
  where $\patch K$ denotes the patch of $K$. Using a discrete
  Cauchy-Schwarz inequality
  \begin{equation}
    \begin{split}
      \cI_3
      &\leq
      C
      \qp{\sum_{K\in\T{}} \qb{\eta_{R}^2 + \sum_{e\in\partial K}\eta_{J}^2}}^{1/2}
      \qp{\sum_{K\in\T{}} \Norm{\nabla v}_{\leb{2}(\patch K)}^2}^{1/2}
      \\
      &\leq
      C \qp{\sum_{K\in\T{}} \qb{\eta_{R} + \sum_{e\in\partial K} \eta_{J}^2} }^{1/2}
      \Norm{\nabla v}_{\leb{2}(\W)},
    \end{split}
  \end{equation}
  hence, making use of (\ref{eq:first}) and (\ref{eq:second}), we have
  \begin{equation}
    \ltwop{\mathfrak R}{v}
    =
    \cI_1+\cI_2+\cI_3
    \leq
    C \qp{\sum_{K\in\T{}} \qb{\eta_{R} + \sum_{e\in\partial K} \eta_{J}^2} }^{1/2}
    \Norm{\nabla v}_{\leb{2}(\W)},
  \end{equation}
  where the constant $C$ depends only upon the shape regularity of the
  mesh, $p$ and $u_h$. The result follows by dividing through by
  $\Norm{\nabla v}_{\leb{2}(\W)}$ and taking the supremum over all
  possible $0\neq v\in\hoz$.
\end{Proof}

\begin{Cor}
  \label{cor:apost}
  Making use of the triangle inequality, one may show under the
  conditions of Theorem \ref{the:apost-E} the following result holds:
  \begin{equation}
    \enorm{u - u_h}^2
    +
    2C_L \Norm{u-u_h}_{\leb{p}(\W)}^p
    \leq
    C\sum_{K\in\T{}} \qb{ \eta_R^2
    +
    \sum_{e\in \partial K} \eta_J^2}.
  \end{equation}
\end{Cor}

\section{Numerical experiments}
\label{sec:numerics}

We now illustrate the performance of the scheme through a series of
numerical experiments. 

\subsection{Test 1 -- Asymptotic behaviour approximating a smooth solution}
\label{sec:test1}

\begin{figure}[h!]
  \caption[]
  {\label{fig:convergence-smooth}
    Convergence plots for the dG scheme (\ref{eq:dg}) for Test 1. We
    measure error norms involving the dG solution, $u_h$ and the a
    posteriori estimator given in Corollary \ref{cor:apost}.
  }
  \begin{center}
    \subfloat[{
        {
          $k=1, p=4$
        }
    }]{
      \includegraphics[scale=\figscale,width=0.45\figwidth]{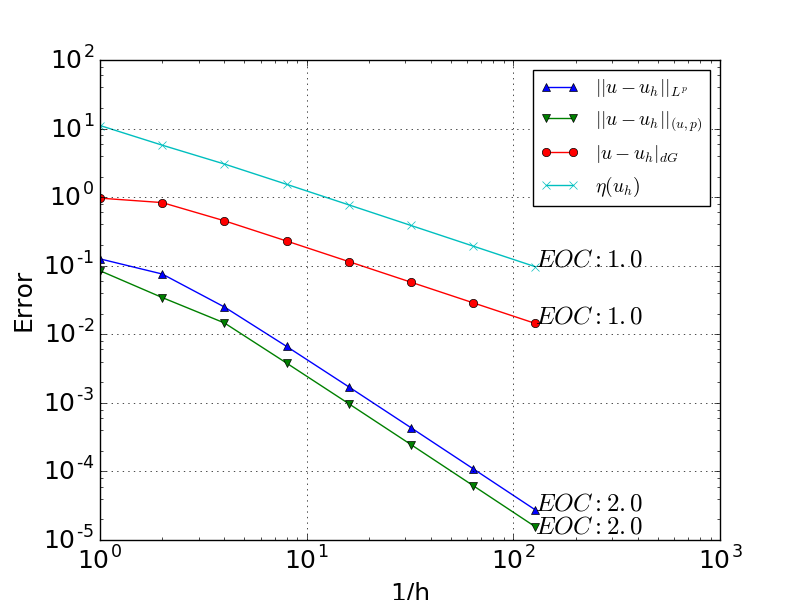}
    }    
    \subfloat[{
        {
          $k=2, p=4$
        }
    }]{
      \includegraphics[scale=\figscale,width=0.45\figwidth]{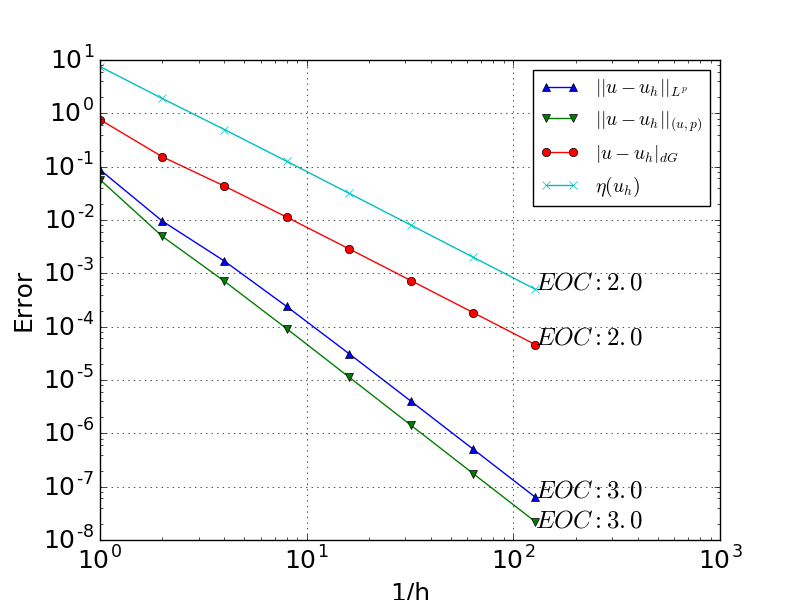}
    }\\
        \subfloat[{
        {
          $k=1, p=8$
        }
    }]{
      \includegraphics[scale=\figscale,width=0.45\figwidth]{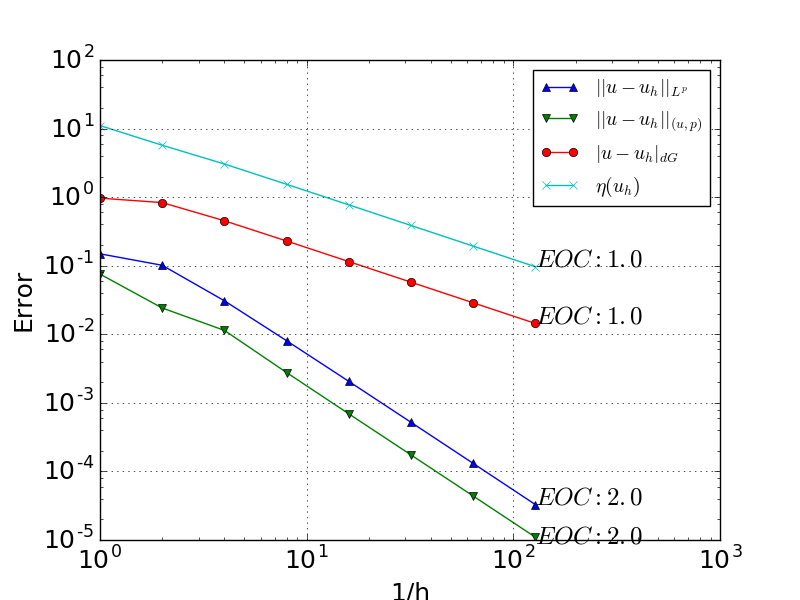}
    }
      \subfloat[{
        {
          $k=2, p=8$
        }
    }]{
      \includegraphics[scale=\figscale,width=0.45\figwidth]{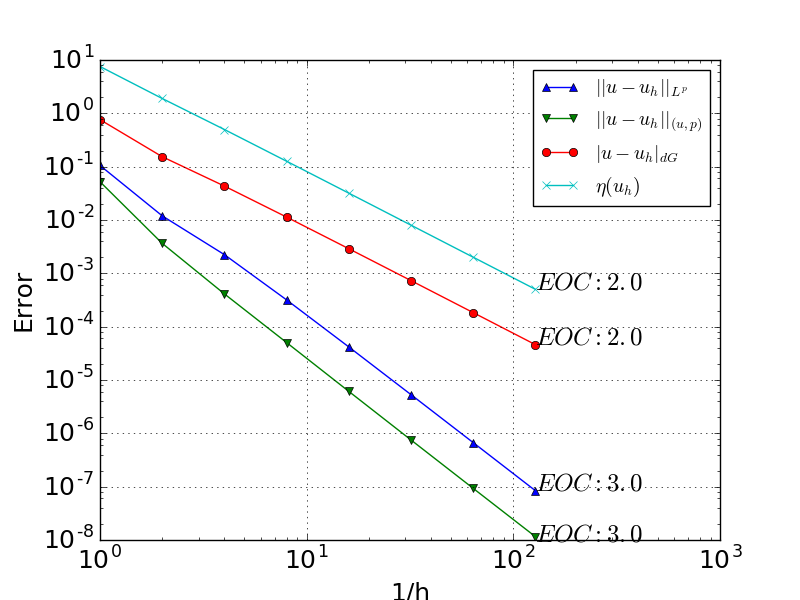}
    }
  \end{center}
\end{figure}

As a first test, we consider the domain $\Omega = [0,1]^2$. We fix $f$
such that the exact solution is given by
\begin{equation}
  \label{eq:smooth-sol}
  u(x,y) = \sin{\pi x}\sin{\pi y},
\end{equation}
and approximate $\Omega$ through a uniformly generated, criss-cross
triangular type mesh to test the asymptotic behaviour of the numerical
approximation. The results are summarised in Figure
\ref{fig:convergence-smooth} (A) -- (D), and confirm the theoretical
findings in Sections \ref{sec:quasinorms} and \ref{sec:apost}.

More specifically, we consider the case $k=1,2$, $p=4,8$ and show that
convergence measured the $\leb{p}$-norm, the $(u,p)$ quasinorm and the
dG norm are all optimal. Notice that the fact the $\leb{p}$ norm is
optimal is contrary to the analysis. This is a well known fact
\cite{Pryer:2018,KatzourakisPryer:2018a}. In addition, the a
posteriori estimator is of optimal rate with an effectivity index of
just under 10.

\subsection{Test 2 -- Behaviour of an adaptive scheme for various values of $p$}
\label{sec:test2}

We consider the domain $\Omega = [0,1]^2$ and fix $f = 1000$ in this
case there is no known solution. However, examining the energy
functional (\ref{eq:energy}) one can see that a minimiser has to
'balance' the $\leb{2}$ norm of its derivative with the $\leb{p}$ norm
of the function. For large $p$ this almost translates into control of
the $\esssup$ which causes boundary layers to appear.

We approximate $\Omega$ through a uniformly generated, criss-cross
triangular type initial mesh consisting of $4$ elements. We run an
adaptive algorithm of {SOLVE}, {ESTIMATE}, {MARK}, {REFINE} type,
where {SOLVE} consists of solving the formulation (\ref{eq:dg}),
{ESTIMATE} is done through the evaluation of the estimator given in
Corollary \ref{cor:apost}, {MARK} is a maximum strategy with $50\%$ of
the elements marked for refinement at each iteration and {REFINE} is a
newest vertex bisection.

\begin{figure}[h!]
  \caption[]
          {\label{fig:adapt}
            Solutions of the dG scheme (\ref{eq:dg}) for Test 2 and the underlying meshes generated through the adaptive algorithm described in the text. Notice that as $p$ increases, boundary layers become present in the solution requiring mesh resolution to be concentrated at the boundary.
          }
  \begin{center}
    \subfloat[{
        {
          $p=2$. The mesh consists of $924,085$ elements after $13$ iterations.
        }
    }]{
      \includegraphics[scale=\figscale,width=0.45\figwidth]{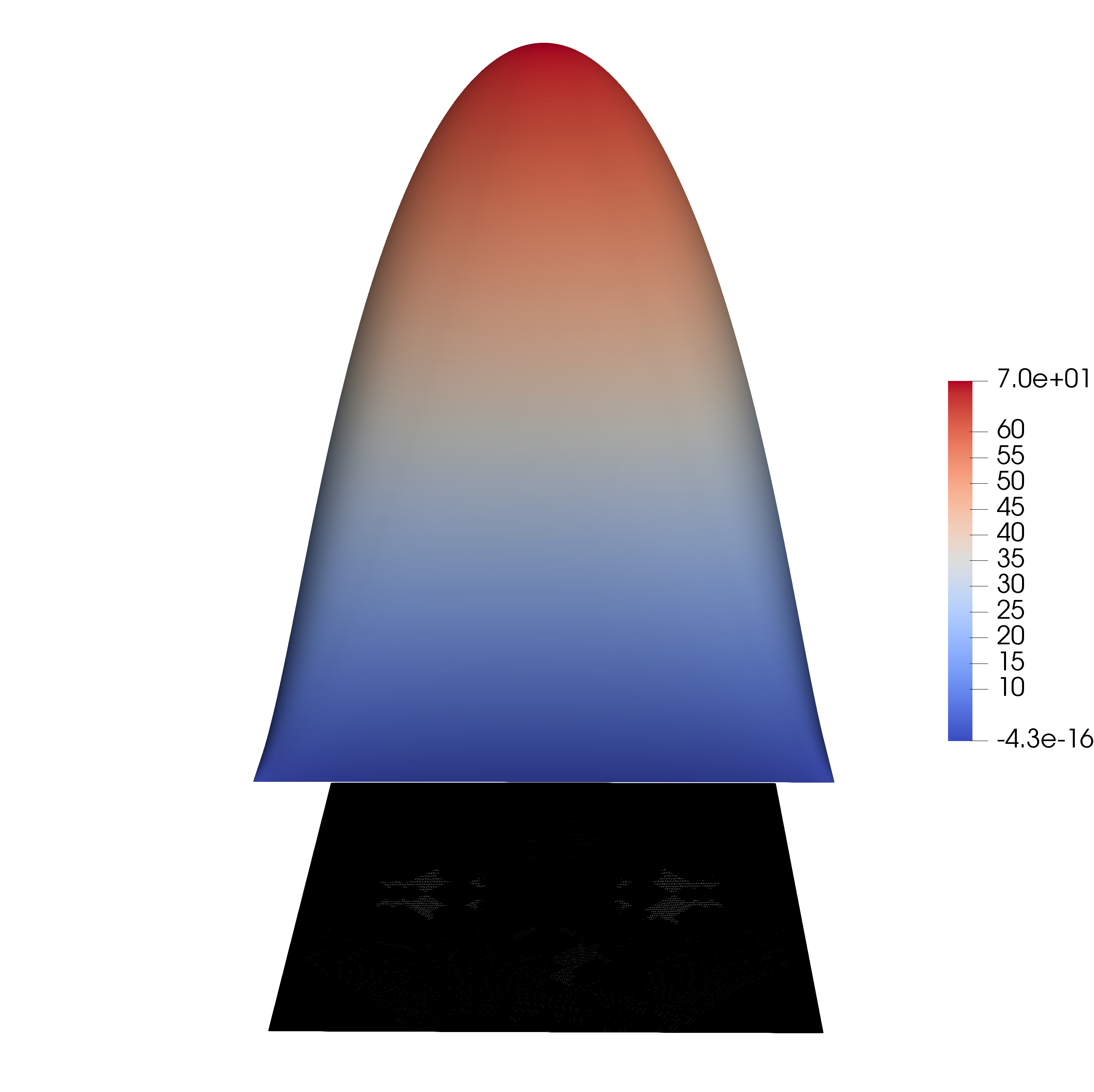}
    }
    \subfloat[{
        {
          $p=4$. The mesh consists of $361,680$ elements after $13$ iterations.
        }
    }]{
      \includegraphics[scale=\figscale,width=0.45\figwidth]{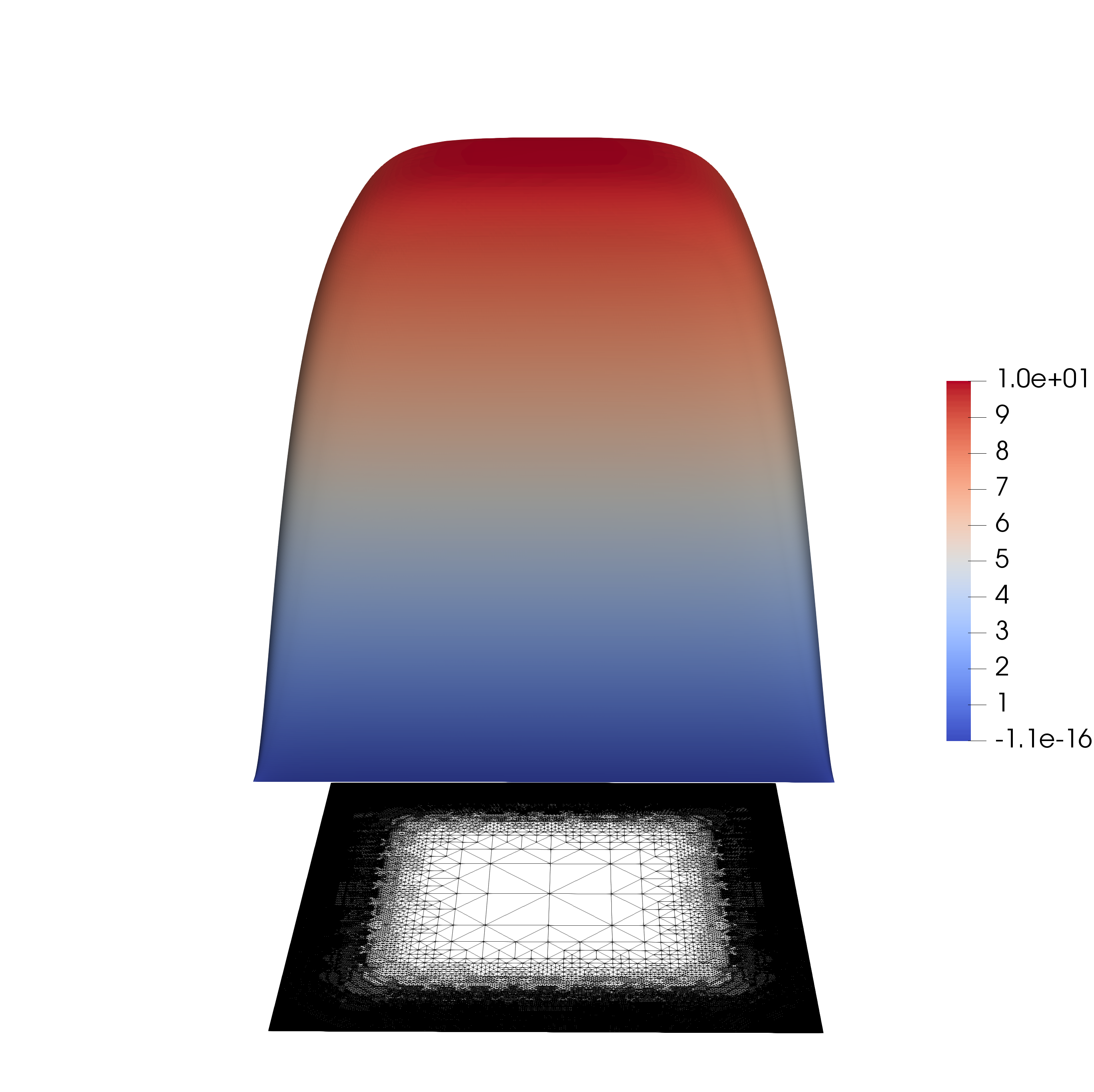}
    }
    \\
    \subfloat[{
        {
          $p=8$. The mesh consists of $339,392$ elements after $13$ iterations.
        }
    }]{
      \includegraphics[scale=\figscale,width=0.45\figwidth]{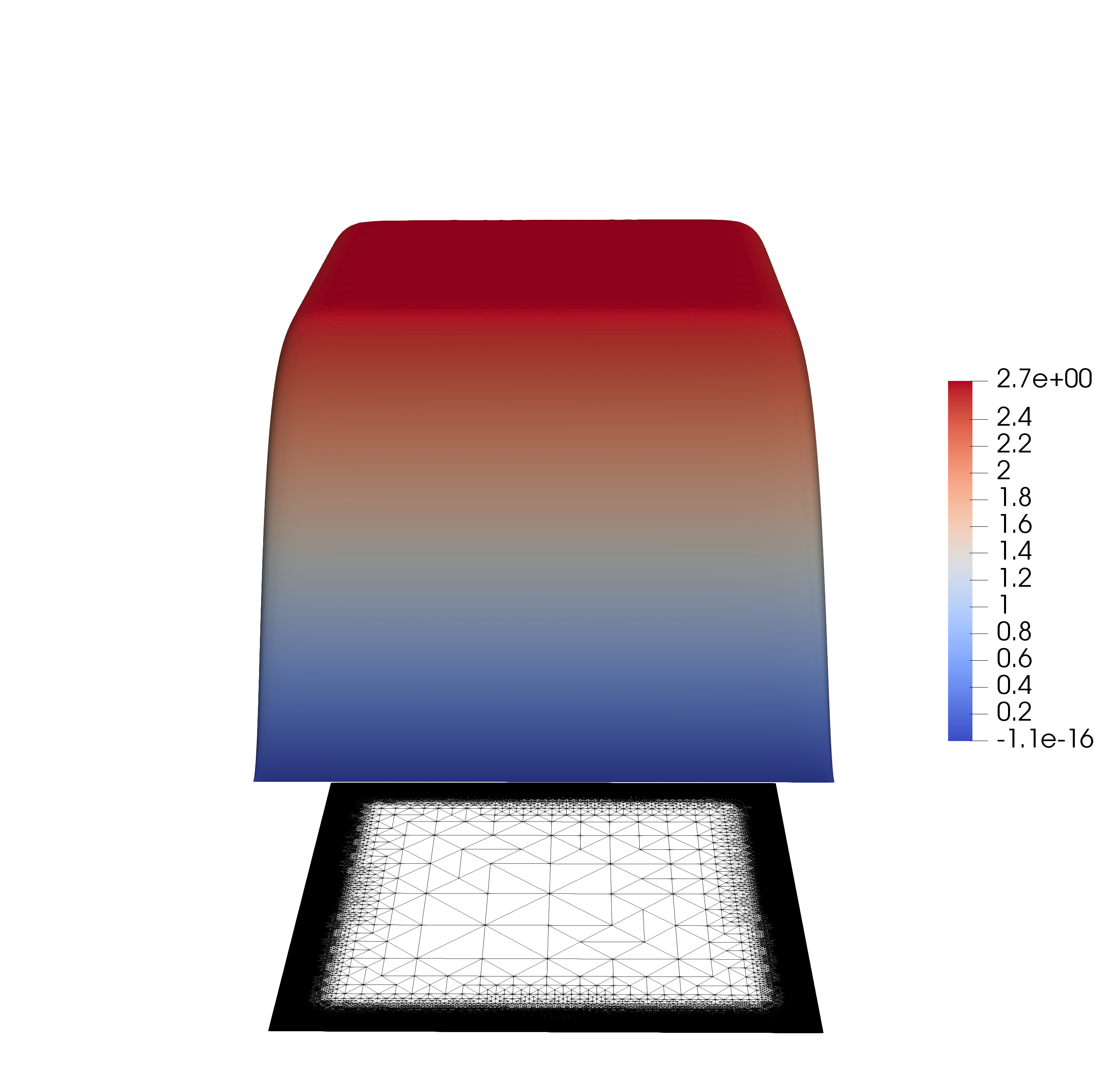}
    }
      \subfloat[{
        {
          $p=12$. The mesh consists of $340,176$ elements after $13$ iterations.
        }
    }]{
      \includegraphics[scale=\figscale,width=0.45\figwidth]{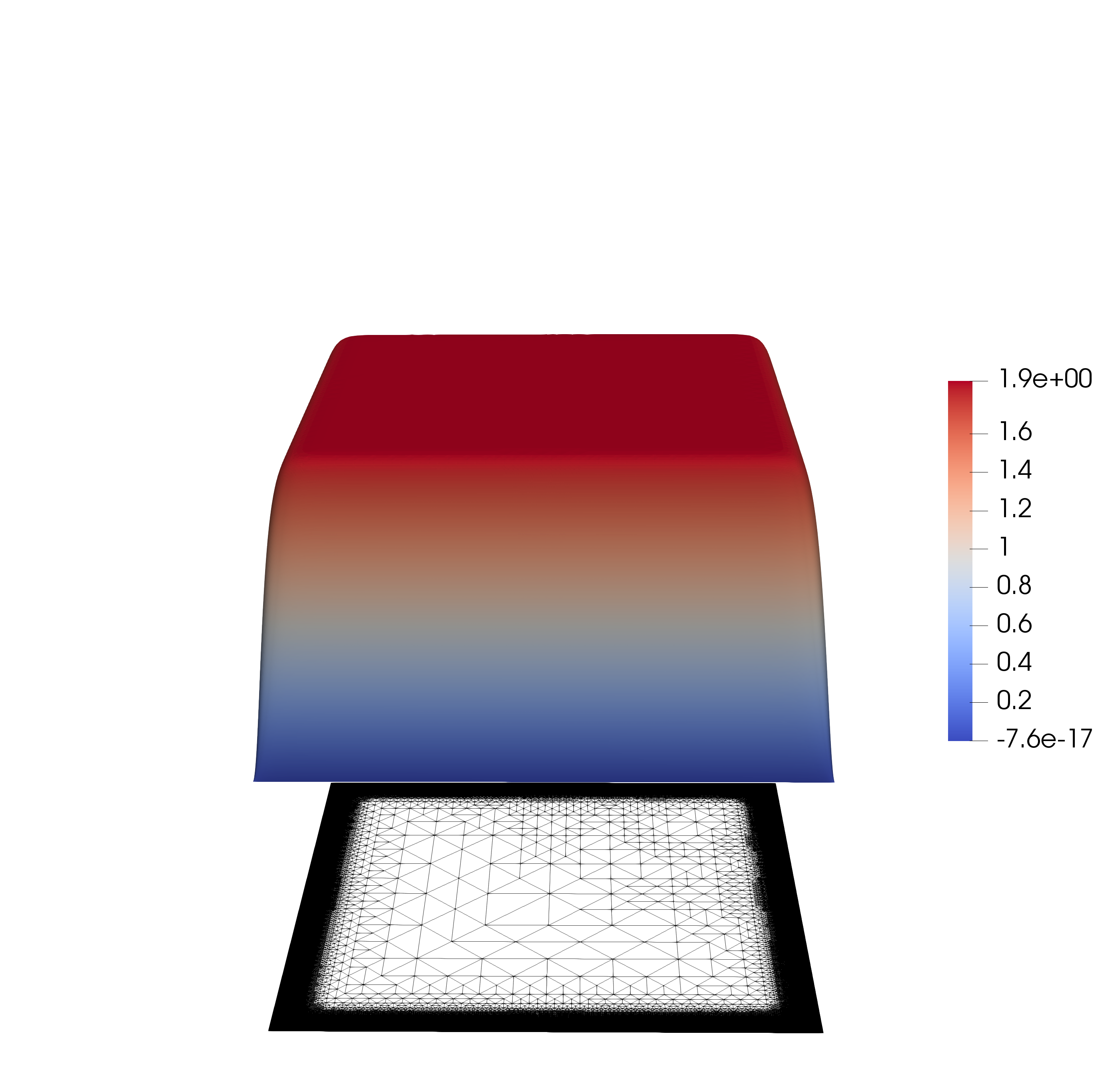}
    }
  \end{center}
\end{figure}

The results are summarised in Figure \ref{fig:adapt} (A) -- (D) where
we consider the case $k=1$, $p=2,4,8,12$ and examine the solution
and underlying adaptive mesh.

\end{document}